\magnification=\magstep1
\input amstex
\documentstyle {amsppt}
\nologo
\pageheight {8.5 truein} \pagewidth {6.0 truein} \NoBlackBoxes
\NoRunningHeads \topmatter
\define\card{\operatorname{card}}

\title
Factorizations of polynomials over noncommutative algebras and
sufficient sets of edges in directed graphs
\endtitle
\author
Israel Gelfand, Sergei Gelfand,  Vladimir Retakh, and Robert Lee
Wilson
\endauthor
\abstract To directed graphs with unique sink and source we
associate a noncommutative associative algebra and a polynomial
over this algebra. Edges of the graph correspond to pseudo-roots
of the polynomial. We give a sufficient condition when
coefficients of the polynomial can be rationally expressed via
elements of  a given set of pseudo-roots (edges). Our results are
based on a new theorem for directed graphs also proved in this
paper.
\endabstract
\email
\newline igelfand\@math.rutgers.edu, sxg\@math.ams.org,
\newline vretakh\@math.rutgers.edu, rwilson\@math.rutgers.edu
\endemail
\address
\newline I.G., V.R., R.W.: Department of Mathematics, Rutgers
University, Piscataway, NJ 08854-8019
\newline S.G.: American Mathematical Society, 201 Charles Street,
Providence, RI 02904-2213
\endaddress
\thanks
Vladimir Retakh was partially supported by NSA
\endthanks
\dedicatory
To the memory of Felix Alexandrovich Berezin
\enddedicatory
\endtopmatter
\document
\bigskip
\head 0. Introduction\endhead
\medskip
Let $R$ be an associative ring with unit and
$P(t)=a_0t^n+a_1t^{n-1}+\dots +a_n$ be polynomial in $R[t]$. Here
$t$ is an independent central variable. Recall \cite {GGRSW} that
an element $x\in R$ is {\it a pseudo-root} of $P(t)$ is there
exist polynomials $Q_1(t), Q_2(t)\in R[t]$ such that
$$P(t)=Q_1(t)(t-x)Q_2(t).$$

The element $x$ is {\it a right root} of $P(t)$ if $Q_2(t)=1$ and
is {\it a left root} of $P(t)$ if $Q_1(t)=1$. It is
easy to check that $x$ is a right root of $P(t)$ if and only if
$a_0x^n+a_1x^{n-1}+\dots +a_n=0$. Similarly, $x$ is a left root if
and only if $x^na_0+x^{n-1}a_1+\dots +a_n=0$.

For noncommutative $R$, a theory of polynomials
over $R$ should be based not only on properties of right (left)
roots but on pseudo-roots as well.

 Suppose now that $P(t)$ is a monic polynomial (i. e., 
$a_0=1$). We say that a subset $Y$ of pseudo-roots of $P(t)$ is
{\it a defining set} if $P(t)$ can  be factored as
$P(t)=(t-y_n)(t-y_{n-1})\dots (t-y_1)$, where $y_k\in Y$ for
$k=1,2,\dots , n$.

When $R$ is a commutative integral domain, 
the set of pseudo-roots of $P(t)$
coincides with the set of  roots $X$. In this case if  $\card (X)=n$ 
then $X$ is a defining set and the coefficients of
$P(t)$ can be written as polynomial expressions in $x_1,
x_2,\dots , x_n\in X$ (the Vi\'ete theorem).

When $R$ is not commutative one faces new phenomenona.
To describe these we introduce the following notation.
For a subset $Y\subseteq R$,
denote by $R(Y)$ the subring generated by $Y$ and by $\tilde
R(Y)$ the subring  of all $r\in R$ such that $r$ can be written
as a rational expression in $y\in Y$. Then

\roster

\item
The coefficients of $P(t)$ can be written only as {\it rational}
expressions in left or right roots of $P(t)$ (see \cite {GR1} and
other papers \cite {GR2, GGRW}). More precisely, $a_1,a_2,\dots ,
a_n$ can be written as polynomial expressions in $z\in Z$ where
$Z$ is a defining set of pseudo-roots of $P(t)$, but elements
$z\in Z$ can only be expressed {\it rationally} in terms of
left or right roots of $P(t)$;

\item
There might be a large number of pseudo-roots: if a polynomial
$P(t)$ of degree $n$ has $n$ right roots in a generic
position (see section 1.2) then $P(t)$ has at least $n2^{n-1}$
pseudo-roots;

\item Not any subset $Y$ of pseudo-roots of $P(t)$ of cardinality
greater or equal $n$ is a defining set.
Moreover, the subring $\tilde R(Y)$ (and,
therefore, $R(Y)$)  may not even contain a defining set.
\endroster

In this paper we study the following problem: describe subsets
$Y$ of pseudo-roots of a monic polynomial $P(t)$ such that the
subring $\tilde R(Y)$ contains a defining subset of pseudo-roots
of $P(t)$. In fact, we replace $\tilde R(Y)$ by a smaller subring
containing $R(Y)$. However, item (1) above indicates that
a consideration of subrings $R(Y)$ for this problem is too
restrictive.

Surprisingly enough, for a large class of noncommutative algebras
the answer may be given in terms of properties of directed graphs.

The content of the paper is the following. In Section 1 we
formulate our results for the universal algebra of pseudo-roots
of polynomials over noncommutative algebras. In Section 2 we
introduce a class of noncommutative algebras associated with
directed graphs and show their relations with factorizations of
noncommutative polynomials. In Section 3 we introduce {\it
sufficient sets} of edges in directed graphs, study properties of
the sufficient sets of edges, and apply those properties for a
description of sufficient sets of pseudo-roots of noncommutative
polynomials. We believe that the results obtained in this section
are of an interest for a ``pure theory" of directed graphs.

\head 1. Main example: Universal algebra of pseudo-roots\endhead

\subhead {1.1. Algebra $Q_n$}\endsubhead The universal algebra of 
pseudo-roots $Q_n$ over a field $k$ was introduced in \cite {GRW}
and studied in \cite {GGR, GGRSW, SW, P}. Algebra $Q_n$ was used to
construct a natural family of defining sets of pseudo-roots
for noncommutative polynomials.

Recall that generators of $Q_n$ are elements $x_{A,i}$ where
$A\subseteq \{1,2,\dots ,n\}$, $i=1,2,\dots , n$ and $i\notin A$.
The defining relations in $Q_n$ are as follows:
$$
x_{A\cup \{i\},j}+x_{A,i}=x_{A\cup \{j\}, i}+x_{A, j}, \tag 1.1
$$
$$
x_{A\cup \{i\},j}\cdot x_{A,i}=x_{A\cup \{j\},i}\cdot x_{A, j} \tag1.2
$$
for each $A\subseteq \{1,2,\dots , n\}$, $i,j\notin A$, $i\neq j$.  We
will call the generators $x_{A,i}$ the pseudo-roots.

There is a canonical polynomial $\Cal P(t)\in Q_n[t]$ having
$x_{A,i}$'s as its pseudo-roots. Let $\{i_1, i_2,\dots , i_n\}$
be any ordering of $\{1,2,\dots , n\}$. Set $A_1=\emptyset $,
$A_2=\{i_1\}$, $A_3=\{i_1,i_2\}$, $\dots $, $A_n=\{i_1, i_2,\dots
, i_{n-1}\}$. Set
$$
\Cal P(t)=(t-x_{A_n,i_n})(t-x_{A_{n-1}, i_{n-1}})\dots (t-x_{A_1,
i_1}).
$$

\proclaim {Theorem 1.1.1} (\cite {GRW}) Polynomial $\Cal P(t)$
does not depend on an ordering of $\{1,2,\dots , n\}$.
\endproclaim

\proclaim {Theorem 1.1.2} (\cite {GRW}) There exists the canonical
epimorphism $$c: Q_n\rightarrow F[t_1, t_2,\dots , t_n]$$ given by
the formulas $c(x_{A,i})=t_i$.
\endproclaim

Here $F$ is the ground field, $t_1, t_2, \dots , t_n$ are
independent commuting variables.

 \subhead {1.2. Universality of algebra $Q_n$}\endsubhead
Let $R$ be an associative algebra and $P(t)\in R[t]$ be a monic
polynomial of degree $n$. Let $X=\{x_1,x_2, \dots , x_n\}$ be a
set of right roots of $P(t)$. We call the set $X$ {\it generic} if
the Vandermonde matrix
$$
V(i_1, i_2, \dots  ,i_{k+1})
=\left ( \matrix x_{i_1}^k&x_{i_2}^k&\dots &x_{i_{k+1}}^k\\
                           &         &\dots &           \\
                  x_{i_1}  &x_{i_2} &\dots  &x_{i_{k+1}}   \\
                     1     &1       &\dots      &1\endmatrix
\right )
$$ 
is invertible in $R$ for each $k=1,2,\dots , n-1$ and
distinct $1\leq i_1, i_2,\dots ,i_{k+1}\leq n$.

In this case one can define the Vandermonde quasideterminants
$v(i_1, i_2, \dots  ,i_{k+1})$ (\cite {GR1, GR2, GGRW} as follows.
Set $r$ to be the row matrix $(x_{i_1}^k, x_{i_2}^k, \dots
x_{i_k}^k)$ and $c$ to be the column matrix
$(x_{i_{k+1}}^{k-1},x_{i_{k+1}}^{k-2},\dots ,1)^T$. Then
$$v(i_1, i_2, \dots  ,i_{k+1})=x_{i_{k+1}}^k - rV(i_1, i_2, \dots
,i_k)^{-1}c.$$

\example {Example} $v(1,2)=x_2-x_1$.
\endexample

The following proposition follows from the general properties
of quasideterminants (\cite {GR1, GR2, GRRW}).

\proclaim {Proposition 1.2.1}
 If $X$ is a generic set then all
quasideterminants $v(i_1, i_2, \dots ,i_{k+1})$ are defined and
invertible in $R$.
\endproclaim

\proclaim {Theorem 1.2.2} \cite {GRW} If $X=\{x_1, x_2, \dots ,
x_n\}$ is a generic set of right roots of a monic polynomial
$P(t)\in R[t]$ of degree $n$ then there exists a canonical
homomorphism $\phi :Q_n\rightarrow R$ such that $\phi
(x_{\emptyset , k})=x_k$ for $k=1,2,\dots , n$.

Moreover, 
$$
\phi (x_{\{i_1,i_2,\dots , i_k\}, i_{k+1}})= v(i_1,
i_2, \dots  ,i_{k+1})x_{i_{k+1}}v(i_1, i_2, \dots
,i_{k+1})^{-1}
$$ 
and the canonical extension of $\phi $ to the
homomorphism of $Q_n[t]$ to $R[t]$ maps $\Cal P(t)$ to $P(t)$.
\endproclaim

By abusing notation we denote elements $\phi (x_{\{i_1,i_2,\dots
, i_k\}, i})\in R$ via $x_{i_1,i_2,\dots , i_k, i}$.

 \proclaim {Corollary 1.2.3} For a generic set
$X$ and an ordering $\{i_1, i_2,\dots , i_n\}$ of $\{1,2,\dots ,
n\}$ set $y_1=x_{i_1}$, $y_2=x_{i_1,
 i_2}$, $\dots $, $y_n=x_{i_1i_2\dots i_{n-1}, i_n}$. Then
 $P(t)=(t-y_n)(t-y_{n-1})\dots (t-y_1)$.
 \endproclaim

Corollary 1.2.3 implies that $Y=\{y_1, y_2,\dots , y_n\}$ is a
defining set of pseudo-roots of $P(t)$.

\example {Example} If $n=2$ then
 $P(t)=(t-x_{i,j})(t-x_i)$,
 where $x_{i,j}=(x_j-x_i)x_j(x_j-x_i)^{-1}$,
$i, j=1,2$ and  $i\neq j$.

If $n=3$ then $P(t)=(t-x_{ij, k})(t-x_{i,j})(t-x_k)$, where
$$x_{ij,k}=(x_{i,k}-x_{i,j})x_{i,k}(x_{i,k}-x_{i,j})^{-1}=
(x_{j,k}-x_{j,i})x_{j,k}(x_{j,k}-x_{j,i})^{-1}.$$ Here
$i,j,k=1,2,3$ and $i,j,k$ are all distinct.
\endexample

Note that elements $y_k$, $k=1,2,\dots , n$ are rational
expressions
 in $x_{i_1}, x_{i_2},\dots , x_{i_k}$ and do not depend on
 ordering of $i_1, i_2,\dots , i_{k-1}$.

\subhead 1.3. Algebras $Q_n$ and graphs \endsubhead For each natural $n$ we  
define the directed graph $\Gamma_n$ as follows. The vertices of
$\Gamma _n$ are the subsets of the set $\{1,2,\dots , n\}$ (including
the empty subset). The edges of $\Gamma _n$ are all directed
pairs of subsets $(A\cup \{i\}, A)$, where $A\subseteq \{1,2,\dots ,
n\}$ and $i\notin A$; we denote such an edge by $(A,i)$. The subset $A\cup\{i\}$ is the tail of
the edge $(A,i)$ and the subset $A$ is the head of $(A,i)$.

 There is a one-to-one correspondence between the generators $x_{A,i}$
of algebra $Q_n$ and the edges $(A,i)$ of the graph $\Gamma _n$. Note
that two edges $(A,i)$ and $(B,j)$ have a common head if $A=B$ and a
common tail if $A\cup \{i\}$=$B\cup \{j\}$. Also, the generators
occuring in relations (1.1) and (1.2) for the algebra $Q_n$ correspond
to the edges of the diamond $\Cal D$ in $\Gamma _n$ with vertices $A$,
$A\cup \{i\}$, $A\cup \{j\}$, $A\cup \{i,j\}$.

\definition{Definition 1.3.1} We sat that the edge $(A,i)$ in $\Cal D$ 
is obtained by the $D$-operation from the ordered pair of edges
$(A\cup \{i\}, j)$ and $(A\cup \{j\}, i)$ and the edge $A\cup \{i\},
j)$ is obtained by the $U$-operation from the ordered pair of edges
$(A,i)$, $(A,j)$.
\enddefinition

Note that $(A,j)$ is obtained by the $D$-operation from the pair
of edges $(A\cup \{j\}, i)$ and $(A\cup \{i\}, j)$ and the edge
$A\cup \{j\}, i)$ is obtained by the $U$-operation from the
ordered pair of edges $(A,j)$, $(A,i)$.

\definition{Definition 1.3.2} We say that a pseudo-root 
$\xi \in Q_n$ is obtained from an ordered pair of pseudo-roots
$x_{A,i}, x_{B,j}$ by the $u$-operation if the edges $(A,i),(B,j)$
have a common head and $(x_{A,i}- x_{B,j})x_{A,i}= \xi (x_{A,i}-
x_{B,j})$.

A pseudo-root $\eta \in Q_n$ is obtained from an ordered pair of
pseudo-roots $x_{A,i}, x_{B,j}$ by the $d$-operation if the edges
$(A,i),(B,j)$ have a common tail and $(x_{A,i}- x_{B,j})\eta=
x_{A,i}(x_{A,i}- x_{B,j})$.
\enddefinition


A connection between $d$- and $u$-operations in $Q_n$ and $D$-
and $U$-operations in $\Gamma _n$ is given by the following
proposition.

\proclaim {Proposition 1.3.3} The element $x_{A\cup \{i\}, j}$ is
obtained by the $u$-operation from the pair $x_{A, i}, x_{A,j}$.

The element $x_{A, i}$ is obtained by the $d$-operation from the pair
$x_{A\cup \{j\}, i}, x_{A\cup \{i\}, j}$.
\endproclaim

\demo {Proof} Look at the diamond $\Cal D$ and use formulas (1.1) and
(1.2).
\enddemo

 \subhead 1.4. Sufficient sets of pseudo-roots in $Q_n$
\endsubhead
Let $Z$ be a subset of pseudo-roots in  $Q_n$.

\definition
{Definition 1.4.1} The set of elements in $Q_n$ that can be
obtained from elements of $Z$ by a successive applications of
$d$- and $u$-operations is called the $du$-envelope of $Z$.
\enddefinition

The following proposition is obvious.

\proclaim {Proposition 1.4.2} Let $f : Q_n\rightarrow D$ be a
homomorphism of $Q_n$ into a division ring $D$. If $r\in Q_n$
belongs to a $du$-envelope of $Z$ then $f(r)$ can be written as a
rational expression in elements  belonging to $f(Z)$.

\endproclaim

\remark {Remark} In fact, $f(r)$ can be obtained from elements
$f(z)$, $z\in Z$ by operations of addition, subtractions,
multiplication, left conjugation $(a,b)\mapsto (a-b)a(a-b)^{-1}$
and right conjugation $(a,b)\mapsto (a-b)^{-1}a(a-b)$.
\endremark

Recall that a subset $Y$ of pseudo-roots of $\Cal P(t)$ is {\it a
defining set} is $\Cal P(t)$ can  be factorized as $\Cal
P(t)=(t-y_n)(t-y_{n-1})\dots (t-y_1)$, where $y_k\in Y$ for
$k=1,2,\dots , n$.

\definition {Definition 1.4.3} A set $Z\subseteq Q_n$ is called sufficient
if the $du$-envelope of $Z$ contains a defining set of pseudo-roots of
$\Cal P(t)$.
\enddefinition

By Proposition 1.4.2, instead of asking when the rational envelope of
a set $W$ of pseudo-roots of a monic polynomial $P(t)\in D[t]$
contains a defining set of pseudo-roots of $P(t)$ we ask a more
restrictive question: when the set $W$ is an image of a sufficient set
of pseudo-roots in $Q_n$?

Any defining set of elements is a sufficient set. Other examples
of sufficient sets in $Q_n$ are given by the following statement.

 \proclaim {Proposition 1.4.4} The sets $\{x_{\emptyset,
k},\,k=1,2,\dots , n\}$ and $\{x_{12\dots \hat k \dots n, k},\,
k=1,2,\dots , n\}$ are sufficient in $Q_n$ for $\Cal P$.
\endproclaim

\demo {Proof} Use successively Proposition 1.3.3.
\enddemo

A necessary condition of a subset in $Q_n$ to be sufficient for
$\Cal P(t)$ is given by the following theorem.

\proclaim {Theorem 1.4.5} If $Z=\{x_{A_1, i_1}, x_{A_2, i_2},
\dots , x_{A_n, i_n}\}$ is a sufficient subset of $Q_n$ then all
$i_1, i_2,\dots , i_n$ are distinct.
\endproclaim

\demo {Proof} Let $R$ be a commutative algebra without zero
divisors, $P(t)$ be a monic polynomial over $R$ and $X={x_1,
x_2,\dots , x_n}$ a generic set of right roots of $P(t)$.
According to Theorem 1.2.2 there exists a homomorphism $\phi :
Q_n \rightarrow R$ such that $\phi (x_{A_k, i_k})=x_{i_k}$ for all
$k=1,2,\dots , n$. Let $\hat Z$ be the $du$-envelope of $Z$.
Clearly, $\phi (Z)=\phi (\hat Z)$. Therefore, $\phi (Z)$ is a
defining set in $R$. It implies that $\phi (Z)=X$ and so all
$i_1, i_2,\dots , i_n$ are distinct.
\enddemo

A set of edges in a directed graph is connected if it is connected
in the associated non-directed graph (see Section 3 below for
details).

 \proclaim {Theorem 1.4.6} Let $Z=\{x_{A_1, i_1},
x_{A_2, i_2}, \dots , x_{A_n, i_n}\}$ be a subset of $Q_n$ such
that all $i_1, i_2,\dots , i_n$ are distinct. If the set of edges
$\{(A_1, i_1), (A_2, i_2), \dots , (A_n, i_n)\}$ in $\Gamma _n$
is connected then the set $Z$ is sufficient.
\endproclaim

 Let $f:Q_n \rightarrow D$ be a
homomorphism of $Q_n$ into a division ring $D$, $\hat f:Q_n[t]
\rightarrow D[t]$ be the induced homomorphism of the polynomial
rings and $P(t)=\hat f(\Cal P(t)$.

\proclaim {Corollary 1.4.7} Let $Z=\{x_{A_1, i_1}, x_{A_2, i_2},
\dots , x_{A_n, i_n}\}$ be a subset of $Q_n$ such that all $i_1,
i_2,\dots , i_n$ are distinct and the set of edges $\{(A_1, i_1),
(A_2, i_2), \dots , (A_n, i_n)\}$ in $\Gamma _n$ is connected.

Then all coefficients of $P(t)\in D[t]$ can be obtained from
elements $f(z)$, $z\in Z$ by operations of addition, subtractions,
multiplication, and left and right conjugation.
\endproclaim

 \example {Examples} 
1. The set $X=\{x_{\emptyset , 1}, x_{\emptyset , 2}, \dots ,
x_{\emptyset , n}\}$ is connected (the corresponding edges have a
common head $\emptyset $). Therefore, $X$ is a sufficient set.

2. For $n=2$ the sufficient sets are $\{x_{\{i\},j}, x_{\emptyset
,i}\}$, $\{x_{\emptyset ,j}, x_{\emptyset ,i}\}$, $\{x_{\{i\},j},
x_{\{j\},i}\}$. The sets $\{x_{\{i\},j}, x_{\emptyset ,j}\}$ are
not sufficient. Here $i,j=1,2$, $i\ne j$.

3. Let $n=3$. The set $\{x_{\{1\},2}, x_{\{2\},1}, x_{\{1\},3}\}$ is
sufficient because  
$$
(x_{\{2\},1}-x_{\{1\},2})x_{\emptyset,1}
=x_{\{1\},2}(x_{\{2\},1}-x_{\{1\},2}),
$$
$$
x_{\{12\},3}(x_{\{1\},3}-x_{\{1\},2})
=(x_{\{1\},3}-x_{\{1\},2})x_{\{1\},3}
$$
and $\{x_{\{12\},2}, x_{\{1\},2}, x_{\emptyset , 1}\}$ is a
defining set of pseudo-roots.

The sets $\{x_{\{1\},2}, x_{\{2\},1}, x_{\emptyset, 3}\}$ and
$\{x_{\{1\},3}, x_{\{1\},2}, x_{\emptyset , 1}\}$ also are
sufficient but not defining sets in $Q_3$.

4. The set $W=\{x_{\{12\},3}, x_{\{3\},2}, x_{\emptyset ,1}\}$ is not
sufficient because $d$- and $u$-operations are not defined on
elements of $W$. We believe that coefficients of $P(t)$ cannot be
written as rational expressions in elements of $W$, but we do not
have a proof of this statement.
\endexample

A proof of Theorem 1.4.6 follows from a more general theorem for
algebras associated with directed graphs (quivers). We will
describe such algebras in the next section.

\head 2. Algebras associated with directed graphs\endhead

The class of algebras defined in this section contains the
universal algebra of pseudo-roots $Q_n$.

\subhead 2.1. Directed graphs\endsubhead In
this subsection we recall some well-known definitions and results
about directed graphs.

A directed graph $\Gamma $ is a pair of sets $\Gamma =(V, E)$ and a
map $\phi : E\rightarrow V\times V$. Elements of $V$ are called {\it
vertices} of $\Gamma $ and elements $e\in E$ are called {\it edges} of
$\Gamma $.

Let $\phi (e)=(t(e), h(e))$. The vertex $t(e)$ is called the {\it
tail} of $e$ and the vertex $h(e)$ is called the {\it head} of $e$.

A directed path $P$ of length $k$ in $\Gamma $ from a vertex $u$ to a
vertex $v$ is a sequence of edges $e_1, e_2,\dots , e_k$ such that
$t(e_{i+1})=h(e_i)$ for $i=1,2,\dots , k-1$ and $t(e_1)=u$,
$h(e_k)=v$. The vertex $u$ is denoted by $t(P)$ and called the {\it
tail} or {\it origin} of $P$. The vertex $v$ is denoted by $h(P)$ and
called the {\it head} or {\it terminus} of $P$. The length of $P$ is
denoted by $l(P)$.

An edge $e\in E$ is called {\it essential} if there is no path $P$ such
that $l(P)\geq 2$ and $t(e)=t(P)$, $h(e)=h(P)$. A directed path
$P$ is called {\it maximal} if all its edges are essential.

A vertex $u\in V$ is called a {\it source} if there is no edge
$e\in E$ such that $h(e)=u$. A vertex $v\in V$ is called a {\it 
sink} if there is no edge $f\in E$ such that $t(f)=v$.

A directed graph $\Gamma =(V, E)$ is called a {\it layered graph}
if there is a function $r$ from $V$ to the set of non-negative
integers $\Bbb Z_+$ such that $r(t(e))-1=r(h(e))$ for any edge
$e\in E$.

\definition {Definition 2.1.1} A directed graph is called a modular graph if:

\roster
\item  For any two edges $e_1, e_2$ with a common tail there exist
edges $f_1, f_2$ with a common head such that $h(e_i)=t(f_i)$ for
$i=1,2$;
\item  For any two edges $h_1, h_2$ with a common head there exist
edges $g_1, g_2$ with a common tail such that $h(g_i)=t(h_i)$ for
$i=1,2$.
\endroster
\enddefinition
We do not require the uniqueness of $f_1, f_2$ and $h_1, h_2$.

\subhead 2.2. Examples of directed graphs\endsubhead To any
partially ordered set $I$ there corresponds the directed graph $\Gamma
_I$ called the {\it Hasse graph} of $I$. Its vertices are elements
$x\in I$ and its edges are pairs $e=(x, y)\in I\times I$ such that $y$
is an immediate predecessor of $x$ (in other words, $y<x$ and there is
no $z\in I$ such that $y<z<x$.) The element $x$ is the tail of $e$ and
the element $y$ is the head of $e$.

Recall that a partially ordered set $I$ with a function
$r:I\rightarrow \Bbb Z_+$ is called a {\it ranked poset} with the {\it
ranking function} if $r(x)>r(y)$ for any $x>y$. Then $r$ turns the
corresponding Hasse graph $\Gamma _I$ into a layered graph.

To any maximal (minimal) element in $I$ there corresponds a source
(sink) in $\Gamma _I$.

Standard examples of ranked partially ordered sets are as follows.

1. To any set $S$ corresponds the partially ordered set $\Cal
P(S)$ of all subsets of $S$. The order relation is given by
inclusion and $r(A)$ equals the cardinality of $A$.

2. A family $\Cal F\subseteq \Cal P(S)$ is called {\it a complex} if
$B\in \Cal F$ and $A\subseteq B$ imply $A\in \Cal F$.  The order and
the ranking function on $\Cal P(S)$ induce an order and a ranking
function on $\Cal F$.

3. To any finite-dimensional vector space $E$ over a field $\Cal K$
corresponds the partially ordered set $\Cal V(E)$ of all vector
subspaces of $E$. The order relation is given by inclusion and $r(V)$
is the dimension of $V$.

4. To any natural number $n$ there corresponds the partially ordered
set of partitions of $n$, i.e., weakly decreasing sets of natural
numbers $\lambda =(\lambda _1,\lambda _2,\dots , \lambda _k)$ such
that $n=\lambda _1+\lambda _2+\dots+ \lambda _k$. The cardinality of
$\lambda $ is called the {\it length} of the partition. Recall that
$\lambda \leq \mu$ if $\mu =(\mu _1,\dots , \mu _{i_1}, \mu
_{i_1+1},\dots , \mu _{i_2},\mu _{i_2+1},\dots , \mu _{i_k})$ and
$\lambda _1=\mu _1+\dots + \mu _{i_1}$, $\lambda _2=\mu _{i_1+1}+\dots
+\mu_{i_2}$, \dots, $\lambda _k=\mu_{i_{k-1}+1}+\dots+\mu_{i_k}$. 
The ranking function $r(\lambda )$ is the length of $\lambda$.

The corresponding Hasse graphs are layered graphs.  Another
important example of a layered graph is  the graph of right
divisors of a monic polynomial.

\subhead 2.3. The graph of right divisors \endsubhead Let $P(t)$
be a monic polynomial over an associative algebra $R$ and $S$ be a
set of pseudo-roots of $P(t)$. Denote by $R_S$ the subalgebra in
$R$ generated by pseudo-roots $x\in S$.

Construct a layered graph $\Gamma (P, S) =(V, E)$,
$$
V=V(n)\sqcup V(n-1)\sqcup \dots V(1)\sqcup V(0)
$$ 
as follows. The vertices of $V(k)=\{v\in V: r(v)=k\}$ are monic polynomials
$B(t)\in R[t]$ such that $\deg B(t)=k$ and
$$
P(t)=Q(t)B(t)
$$ 
in $R[t]$.

We say that there is an edge from vertex $B_1(t)$ to $B_2(t)$ in
$\Gamma $ if 
$$
B_1(t)=(t-x)B_2(t)
$$ 
for some $x\in S$.

Note that $V(n)$ consists of one vertex $v=P(t)$ and $V(0)$ consists
of one vertex $w=1$.

\subhead 2.4. From graphs to algebras and polynomials\endsubhead
Let $\Gamma =(V, E)$ be a directed graph. Fix a field $k$. Let
$T(E)$ be the free associative algebra over $k$ with generators
$e\in E$. To any directed path $\pi = (e_1, e_2,\dots , e_k)$ in
$\Gamma $ there corresponds a polynomial $U_{\pi }(t)\in T(E)[t]$
$$
U_{\pi }(t)=(t-e_1)(t-e_2)\dots (t-e_k).
$$

The following definition was introduced in \cite {GRSW}.

\definition {Definition 2.4.1}
The algebra $A(\Gamma )$ is the quotient algebra of $T(E)$ modulo
the following relations:
$$ 
U_{\pi _1}(t)=U_{\pi _2}(t) \tag 2.4.1 
$$
in $T(E)$ for any two paths $\pi _1$ and $\pi _2$ that common
beginning and common end.
\enddefinition

For two vertices $v, w$ in $\Gamma $ define the following polynomial
$\Cal P_{v, w}(t)$ over $A(\Gamma )$: 

If there exists a path $\pi =(f_1, f_2,\dots , f_m)$ from $v$ to $w$,
set
$$
\Cal P_{v,w}(t)=(t-f_1)(t-f_2)\dots (t-f_m)
$$
in $A(\Gamma )[t]$. If there are no paths from $v$ to $w$, we set $\Cal
P_{v,w}=1$.

The construction of algebra $A(\Gamma )$ implies that the polynomial
$\Cal P_{v, w}(t)$ does not depend of the choice of a path from $v$ to
$w$.

If the set of vertices of $\Gamma $ contains one source $u$ and one
sink $z$, then the polynomial $\Cal P(t)=\Cal P_{u,z}(t)$ is called
the polynomial associated to graph $\Gamma $.

The following theorem implies that in many important cases
algebras $A(\Gamma )$ are defined by linear and quadratic
relations. These relations correspond to diamonds in $\Gamma $.

\proclaim {Theorem 2.4.2} If $\Gamma $ is a modular layered graph
then defining relations for algebra $A(\Gamma )$ can be given by
formulas (2.4.1) when $l(P_1)=l(P_2)$.
\endproclaim

\example {Example} If $\Gamma $ is the Hasse graph for the set of
all subsets of $\{1,2,\dots , n\}$ then $A(\Gamma )$ coincides
with the universal algebra of pseudo-roots $Q_n$.
\endexample

\remark {Remark} It is natural to study quotient algebras
$A_0(\Gamma )$ of algebras $A(\Gamma )$ given by relations $ef=0$
for all pairs of edges $e,f$ in $\Gamma$ such that $h(e)\neq t(f)$.
\endremark

\subhead {2.5. Universality of $A(\Gamma )$}\endsubhead 

Let $R$ be
an algebra, $P(t)$ a monic polynomial of degree $n$ over $R$, and
$S$ a set of pseudo-roots of $P(t)$. Let $\Gamma (P,S)$ be the
layered graph constructed in Section 2.3. It contains two marked
vertices: the source $v=P(t)$ and the sink $w=1$.

Assume that the set $S$ contains pseudo-roots $a_1, a_2,\dots ,
a_n$ such that $$P(t)=(t-a_1)(t-a_2)\dots (t-a_n).$$ Then the
graph $\Gamma (P,S)$ contains a directed path from $v=P(t)$ to
$w=1$.

Following Section 2.4, we we construct algebra $A(\Gamma (P,S))$ and
polynomial $\Cal P(t)$ over this algebra.

\proclaim {Theorem 2.5.1} There is a canonical homomorphism
$$
\alpha : A(\Gamma (P, S))\rightarrow R \tag 2.5.1 
$$
such that the induced homomorphism of polynomial algebras
$$
\hat \alpha : A(\Gamma (P,S))[t]\rightarrow R[t]
$$ 
maps $\Cal P(t)$ to $P(t)$.
\endproclaim

\demo{Proof} To an edge $e$ in $\Gamma $ there corresponds a pair of
polynomials $B_1(t), B_2(t)$ in $R[t]$ such that $B_1(t), B_2(t)$
divide $P(t)$ from the right and $B_1(t)=(t-a)B_2(t)$.

Set
$$
\alpha (e)=a.
$$
One can see that this map can be uniquely extended to the
homomorphism (2.5.1) and that $\hat \alpha (\Cal P(t))=P(t)$.
\enddemo

\subhead 2.6. Iterated constructions\endsubhead
Let $\Gamma =(V, E)$ be a layered graph, $V=\cup _{k=0}^nV(k)$.
Assume that $V(n)$ contains exactly one vertex $u$, $V(0)$
contains exactly one vertex $z$ and there is a directed path from
$u$ to $z$. Denote by $E_0$ the subset of all edges $e\in E$ such
that there is a directed path from $u$ to $z$ containing $e$.
Denote by $\Gamma _0$ the subgraph of $\Gamma $ generated by
$E_0$.

Following Section 2.4, construct the algebra $A(\Gamma )$ and the
polynomial $\Cal P(t)$. It is clear that $E_0$ can be canonically
identified with a set of pseudo-roots of $\Cal P(t)$. This
identification establishes an isomorphism between the graphs
$\Gamma_0$ and $\Gamma (\Cal P(t))$.

\head 3. Sufficient sets of edges for directed graphs
\endhead

In this section we will define and study sufficient set of edges in
directed graphs $\Gamma =(V,E)$. These sets will provide us with a
construction of sufficient sets of pseudo-roots of polynomials $\Cal
P(t)$ over algebras $A(\Gamma )$. All graphs considered in this
sections are {\it simple} (i.e., if $t(e)=t(f)$ and $h(e)=h(f)$ then
$e=f$) and {\it acyclic} (i.e., there are no directed paths $P$ such
that $t(P)=h(P)$.

\subhead 3.1. $D$- and $U$-operation for directed graphs
\endsubhead
Let $\Gamma =(V, E)$ be a directed graph.

\definition {Definition 3.1.1}
\roster
\item A pair of edges  $f_1, f_2$ with a common head
is obtained from the pair  $e_1, e_2$ with a common tail
by $D$-operation if $h(e_i)=t(f_i)$ for $i=1,2$;

\item A pair of edges  $e_1', e_2'$ with a common tail
is obtained from the pair $f_1', f_2'$ with a common head by
$U$-operation if $h(e_i')=t(f_i')$ for $i=1,2$.
\endroster
\enddefinition

\remark {Remark} We do not require the uniqueness of  $D$- and
$U$-operations.
\endremark

\definition{Definition 3.1.2} A subset $E_0\subseteq E$ is called 
$DU$-complete (or simply complete) if the results of any $D$-operation
or any $U$-operation applied to edges from $E_0$ belong to $E_0$.
\enddefinition

\proclaim {Proposition 3.1.3} For any subset $F\subseteq E$ there
exists a {\it minimal} complete set $\hat F\subseteq E$
containing $F$.
\endproclaim

We call $\hat F$ the {\it completion} of $F$.

To any directed graph $\Gamma =(V, E)$ there corresponds the {\it
double directed graph} $\tilde \Gamma =(V, E\sqcup E_-)$ with the same
set of vertices and the doubled set of edges $E\sqcup E_-$.  Here
$E_-$ is a copy of $E$ and to any edge $e\in E$ there corresponds the
{\it opposite} edge $-e\in E_-$ such that $h(-e)=t(e)$ and
$t(-e)=h(e)$.

Edges $e\in E$ are called {\it positive edges} and edges $-e\in E_-$ are
called {\it negative edges}.

Recall that a path from $v\in V$ to $w\in V$ is a set of edges
$(f_1, f_2,\dots , f_k)$ in the double directed graph $\tilde
\Gamma $ such that $t(f_1)=v$, $h(f_k)=w$ and $t(f_{i+1})=h(f_i)$
for $i=1, 2, \dots , k-1$. A path is called {\it positive} if all
its edges are positive.

Any subset of edges $G$ generates a subgraph of $\Gamma (G)=(V(G), G)$
where $V(G)$ is the union of all heads and tails of edges $e\in G$.

A set of vertices $W\subseteq V$ is {\it connected} if for any vertices $v, w\in
W$ there exists a path from $v$ to $w$. A set of edges $G$ is {\it connected} if
$V(G)$ is connected.

\subhead  3.2. Ample sets and sufficient sets
\endsubhead
 Let $\Gamma =(V, E)$ be a directed graph.

\definition { Definition 3.2.1} A set of edges $G$ in $\Gamma $ is called
{\it sufficient\/} if its completion $\hat G$ contains a positive path
from a source to a sink.
\enddefinition

\definition
{Definition 3.2.2} A set of vertices $W\subseteq V$ is called
{\it ample} if 
\roster
\item For any non-sink vertex $v\in V$ there exists a vertex
$u\in W$ such that
there is no positive path in $\Gamma $ from $u$ to $v$;
\item
For any non-source vertex $v\in V$ there exists a vertex
$w\in W$ such that that there is no positive path in $\Gamma $ from $v$ to $w$.
\endroster
\enddefinition

A set of edges is called ample if the set of its tails and heads is ample.

\smallskip

As an example, consider the graph $\Gamma_n$ of all subsets of $\{1,\dots,n\}$. 
It has one sorce $\{1,\dots,n\}$ and one sink $\emptyset$. 

It is clear that a family 
of vertices $F=\{B_1,\dots,B_K\}$ is ample if and only if for 
each $l\in \{1,\dots,n\}$ there is $B_k\in F$ such 
that $l\in B_k$, and $B_{k'}\in F$ such that $l\notin B_{k'}$.
This immediately provides examples of ample sets of edges in $\Gamma_n$.

\proclaim {Proposition 3.2.3}  A set of edges  $(A_1, i_1),
(A_2, i_2), \dots , (A_n, i_n)$ in $\Gamma_n $ is ample if all $i_1,
i_2,\dots , i_n$ are distinct.
\endproclaim

\demo{Proof} Since all $i_1,\dots,i_n$ are distinct, we can renumber the edges in 
our set so that this set becomes $\{(A_1,1),\dots,(A_n,n)\}$ with $i\notin A_i$. 
Then the corresponding set of vertices is 
$$
\{A_1, A_1\cup \{1\},A_2, A_2\cup \{2\},\dots,A_n, A_n\cup \{n\}\}
$$
(possibly with repetitions), and it is clear that for each $l\in \{1,\dots,n\}$ we have
$$
l\in A_l\cup \{l\},\quad l\notin A_l.
$$
\enddemo

\proclaim {Theorem 3.2.4} Any ample connected set of edges  of a
finite modular directed graph is a sufficient set.
\endproclaim

Our proof of the theorem is based on  the following lemma. Let
$F$ be a complete connected set of edges in a modular directed
graph $\Gamma  =(V, E)$.

\proclaim {Lemma 3.2.5} Let $v, u\in V(F)$ and there is no
positive paths in $\Gamma (F)$ from $u$ to $v$. Then 
\roster
\item There exists an edge $f\in F$ such that $t(f)=v$, 
and
\item
There exists an edge $e\in F$ such that $h(e)=u$. 
\endroster
\endproclaim

\demo{Proof} We will prove the first part of the lemma. The
second part can be proved in a similar way.

Let $\Cal P$ be the set of all shortest paths in $\Gamma (F)$ from
$u$ to $v$.  By assumption, $\Cal P$ does not contain a positive path.
Then any path $P\in \Cal P$ is defined by the sequence of edges
$(\dots, -f_{n-k(P)}, f_{n-k(P)+1},\dots , f_n)$ from $F$.

Take $P\in \Cal P$ such that $k(P)$ is minimal. We claim that
$k(P)=0$. In fact, if $k(P)\neq 0$, then one can take
$$D(f_{n-k(P)+1}, f_{n-k(P)})=(g_{n-k(P)+1}, g_{n-k(P)}).$$

Note that $g_{n-k(P)+1}, g_{n-k(P)}\in F$ because $F$ is complete.
Also, $g_{n-k(P)+1}\neq f_{n-k(P)+2}$ otherwise $P$ is not the
shortest path.

Construct a new path  $P'\in \Cal P$ by replacing the edges $-f_{n-k(P)}, f_{n-k(P)+1}$ with the edges $g_{n-k(P)}, -g_{n-k(P)+1}$.

Then $k(P')=k(P)-1$ and we got a contradiction.

Therefore, $k(P)=0$ and we may set $f=f_n$. This proves the lemma.
\enddemo

\demo {Proof of Theorem 3.2.4} Let $H$ be an ample set of edges in
a finite direct modular graph $\Gamma =(V,E)$ and let $\hat H$ be
the completion of $H$.  Let $P$ be a path of maximal length in
$\Gamma (\hat H)$. Denote by $u$ the beginning of this path and
by $v$ the end of this path.

We claim that $u$ is a source and $v$ is a sink in $\Gamma $.

Indeed,  suppose $v$ is not a sink. Since $G$ is an
ample set, there exist a vertex $w\in V(H)$ such that there is
no positive path in $\Gamma (\hat H)$ from $w$ to $v$. Lemma 3.2.5
implies that there is an edge $f\in \hat H$ such that $t(f)=v$.
By adding $f$ to $P$ one gets a direct path in $\Gamma (\hat H)$
from $u$ to $h(v)$. Then $P$ is not a path of maximal length in
$\Gamma (\hat H)$ and we got a contradiction. It proves that $v$
is a sink.

Similarly, one can prove that $u$ is a source. The theorem is proved.
\enddemo

\subhead {3.3. Sufficient sets of edges and sufficient sets of
pseudo-roots}\endsubhead

Now let $\Gamma =(V, E)$ be a directed graph such that
\roster
\item  $\Gamma $ contains a unique source $M$ and a unique
sink $m$;
\item
 For each vertex $v\in V$ there exist a positive path from $M$ to $v$ and a positive
path from $v$ to $M$.
\endroster
\smallskip
Recall that we associate to $\Gamma $ an algebra $A(\Gamma )$ and
a polynomial $\Cal P(t)\in A(\Gamma )[t]$. The polynomial $\Cal P(t)$ is constructed using a 
positive path from $M$ to $m$ in $\Gamma$, but it does not
depend on the path. To any edge $e\in E$ there corresponds to a
pseudo-root $e\in A(\Gamma )$ of $\Cal P(t)$, and to any positive path
$(e_1, e_2,\dots , e_n)$ from $M$ to $m$ in $\Gamma$ there corresponds the
factorization
$$P(t)=(t-e_1)(t-e_2)\dots (t-e_n) 
\tag 3.3.1
$$
of $\Cal P(t)$ over $A(\Gamma )$.

Let  $e_1, e_2\in E$ be edges with the common tail and
$f_1, f_2\in E$ be edges with the common head such that $h(e_i)=t(f_i)$
for $i=1,2$. From the definition of algebra $A(\Gamma )$ it follows that
$e_1+f_1=e_2+f_2$ and $e_1f_1=e_2f_2$. These formulas imply the following proposition.

\proclaim {Proposition 3.3.1} For $i,j=1,2$, $i\neq j$ we have
$$ e_i(f_i-f_j)=(f_i-f_j)f_i,$$
$$ e_i(e_i-e_j)=(e_i-e_j)f_i.$$
\endproclaim

\proclaim {Corollary 3.3.2} If the pair $f_1, f_2$ is obtained from
the pair by $e_1, e_2$ by $D$-operation then the elements $f_1, f_2$ in
$A(\Gamma )$ are obtained from the elements $e_1, e_2$ by
$d$-operations. If the pair $g_1, g_2$ is obtained from the pair
by $h_1, h_2$ by $U$-operation then the elements $g_1, g_2$ in
$A(\Gamma )$ are obtained from the elements $h_1, h_2$ by
$d$-operations.
\endproclaim

\proclaim {Corollary 3.3.3} If $S\subseteq E$ is an ample set then
there exists a factorization (3.1) of $P(t)$ such that
$du$-completion of $S$ contains elements $e_1, e_2,\dots , e_n$
and, therefore, coefficients of $P(t)$.
\endproclaim

Corollary 3.3.3 and Theorem 3.2.4 implies the following theorem.

\proclaim {Theorem 3.3.4} Let $S\subseteq E$ be an ample connected
set of edges in a modular directed graph $\Gamma =(V,E)$. Then
there exists a factorization (3.1) of $\Cal P(t)$ such that
$du$-completion of $S$ contains elements $e_1, e_2,\dots , e_n$
and, therefore, coefficients of $\Cal P(t)$.
\endproclaim

Let $\Gamma $ be the graph of all subsets of $\{1, 2,\dots , n\}$.
Theorem 1.4.6 follows from Theorem 3.3.4 and Proposition 3.2.3
applied to the graph $\Gamma $.

We will also construct two ``perpendicular" examples of connected
ample sets of edges in $\Gamma $.

\example {Example 1} Take $S=\{x_{\emptyset , 1}, x_{\emptyset ,
2}, \dots , x_{\emptyset , n}\}$. The set  $S$ is connected and
ample. It is easy to see that any edge in $\Gamma $ belongs to
the $DU$-completion of $S$. Therefore, the $du$-completion of $S$
contains elements $e_1, e_2,\dots , e_n$ for any factorization of
$P(t)$.
\endexample
\example {Example 2} Take $T=\{x_{12\dots n-1, n}, x_{12\dots
n-2, n-1},\dots , x_{\emptyset , 1}\}$. The set  $T$ is also
connected and ample. Its $DU$-completion coincides with $T$ and
its $du$-completion defines a factorization of $\Cal P(t)$.
\endexample
Other examples of sufficient sets of pseudo-roots were considered
in Section 1.

\enddocument